\documentclass[12pt,leqno,a4paper]{amsart}
\usepackage{amsmath,amssymb,enumerate, xcolor}

\usepackage{enumerate}

\newtheorem{theorem}{Theorem}[section]
\newtheorem{lemma}[theorem]{Lemma}

\newtheorem*{thmA}{Theorem A}
\newtheorem*{thmB}{Theorem B}
\newtheorem*{thmC}{Theorem C}
\newtheorem*{thmD}{Theorem D}

\theoremstyle{remark}
\newtheorem{remark}[theorem]{Remark}
\newtheorem{example}[theorem]{Example}
\numberwithin{equation}{section}

\begin{document}
\title[ conjugacy classes of  size $2p^n$]{Finite groups with conjugacy classes of size $2p^n$}
     \author[Beltr\'an and Shao]{Antonio Beltr\'an\\
     Departamento de Matem\'aticas\\
      Universitat Jaume I \\
     12071 Castell\'on\\
      Spain\\ 
      \\Changguo Shao\\
      College of science\\
NANJING UNIVERSITY OF POSTS AND TELECOMMUNICATIONS\\
NANJING 210023 YADONG
        \\
      China}

 \thanks{Antonio Beltr\'an: abeltran@uji.es,
 Changguo Shao: shaoguozi@163.com}

\keywords{Burnside $p^a$-Theorem,  conjugacy class sizes, simple groups}

\subjclass[2020]{20D05, 20E45, 20D45}

\begin{abstract}
We extend Burnside's $p^a$-Theorem on conjugacy classes of prime-power size by proving that no non-abelian finite simple group possesses a conjugacy class of size $2p^n$ for any prime $p$ and  $n\geq 0$. Based on this, we partially extend a well-known result of L. Kazarin  by showing  that, given a finite group $G$,   the subgroup generated by a conjugacy class $x^G$ of size $2p^n$  where $x$ has odd order, is a solvable normal subgroup of $G$. If instead the  subgroup  $\langle x^G \rangle$  generated by such a class $x^G$   is non-solvable (which forces $x$ to have even order), then $p= 5$  and  $\langle x^G \rangle/S(\langle x^G \rangle)$ is isomorphic to the symmetric group $S_5$, where $S(\langle x^G \rangle)$ denotes the solvable radical of $\langle x^G \rangle$.
\end{abstract}

\maketitle

\section{Introduction}

Burnside's $p^a$-Theorem establishes that a finite simple group  contains no nontrivial conjugacy class of prime-power size \cite[Theorem 15.2]{Huppert}. To date, no character-free proof  of this fundamental result is known. In 1981, L. Kazarin significantly generalized this  result by showing that every conjugacy class of prime-power size in a finite group $G$   generates a  solvable normal subgroup in $G$,  relying on the weak orthogonality relation for blocks in  Modular Representation Theory (see \cite[Theorem 1]{Kaz} or \cite[Theorem 15.7]{Huppert}). 

\medskip

Building upon these foundational results, in this paper, we extend these results to conjugacy classes of size $2p^n$, where $p$ is a prime and $n$ is  a nonnegative integer. Our first result is the following.

\begin{thmA} \label{A}  Let $G$ be a non-abelian finite simple group. Then $G$ has no conjugacy class of size $2p^n$  for any prime $p$ and any $n\geq 0$.
\end{thmA}

It is worth noting that conjugacy classes with sizes divisible by exactly two primes  do exist among finite simple groups, although rather infrequently. For instance, PSL$(4,3)$ possesses conjugacy  classes of size $2^3\cdot 3^6$ and $2^4\cdot 3$,  while $A_5$ features sizes $2^2\cdot 5$ and $2^2\cdot 3$.  Furthermore, class sizes of the form $2m$, where $m$ is a square-free odd  composite number, also occur in simple groups, such as   $2\cdot 5\cdot 11$ in PSL$(2,11)$,   $2\cdot 3\cdot 7$ in ${\rm PSL}(2,7)$ or  $2\cdot 3\cdot 5\cdot 7$ in $A_8$.

\medskip
In view of Theorem A, we obtain a partial extension of Kazarin's theorem as follows.

\begin{thmB} \label{B} Let $G$ be a finite group and suppose that $1\neq g\in G$ is an element of odd order with $|g^G|=2p^n$, where $p$ is a prime and $n\geq 0$. Then $g^G$ generates a solvable normal subgroup of $G$.
\end{thmB}

The odd-order hypothesis in Theorem B cannot be omitted. Specifically, if $g$ is a transposition of the symmetric group $G= S_5$, then $|g^G|=10$, while $ \langle g^G\rangle =G$ is not solvable. 
 This example allows the construction of further non-solvable groups illustrating that the hypothesis cannot be relaxed (see Section 4). In fact, $A_5$ serves as the core obstruction in all such cases. To formalize this, we first establish:

 \begin{thmC}
     Let $G$ be a non-abelian finite simple group and suppose that there exists $a\in $ {\rm Aut}$(G)$ such that $|G:C_G(a)|=2p^n$ with $p$ prime and $n\geq 1$. Then $G\cong A_5$ , $p=5 $, $n=1$,  and $a$ induces the unique outer automorphism of $A_5$ (corresponding to a transposition of $S_5$).   
    \end{thmC}

As a consequence, we derive the following result.

 \begin{thmD}
 Let $G$ be a finite group containing a conjugacy class $g^G$ of size $2p^n$, where $p$ is prime and $n\geq 1$,  such that $\langle g^G\rangle$ is non-solvable. Then $p=5$, and 
$\langle g^G \rangle/S(\langle g^G \rangle)$ is isomorphic to $S_5$, where $S(\langle g^G \rangle)$ denotes the solvable radical of $\langle g^G \rangle$.
  \end{thmD}

\medskip
The proofs of our results appeal to several deep classification tools, including Character Theory and the Classification of Finite Simple Groups (CFSG).
It remains an open question whether a CFSG-free proof of any of these results exists, or if they depend inherently on the representation theory of simple groups.

\section{Proof of Theorem A}

The proof of Theorem A is divided into three parts according to the three families of simple groups in the CFSG (Theorems \ref{T1}, \ref{T2} and \ref{T3}),  reflecting the fundamentally different nature of the arguments
required in each case. For alternating groups, the proof  only needs elementary  arithmetical computations and Burnside's $p^a$-Theorem. 
\begin{theorem} \label{T1}
The alternating group $A_n$ with $n\geq 5$ does not possess any conjugacy class of size $2p^a$ for any prime $p$ and any integer $a\geq 0$.
\end{theorem}

\begin{proof} 
We can assume $p\neq 2$, as Burnside's $p^a$-Theorem covers the case $p=2$. Let $G=A_n$ and suppose that $K=x^G$ is a conjugacy class of $G$ of size $2p^a$ for some $a\geq 1$. First, we reduce the problem to elements of prime order. Suppose $m = o(x)$ is composite and let $q$ be a prime divisor of $m$. Since $C_G(x) \leq C_G(x^{m/q})$, it follows that $|(x^{m/q})^G|$ divides $|x^G| = 2p^a$. Thus,  again appealing to Burnside's $p^a$-Theorem, we obtain $|(x^{m/q})^G| = 2p^b$ for some $1 \leq b \leq a$. Hence, replacing $x$ by $x^{m/q}$, which has order $q$, we may assume without loss of generality that $o(x)=q$ for some prime $q$.

Therefore, $x$ is a product of $t \geq 1$ disjoint $q$-cycles, where $t$ must be even when $q=2$, and $tq \leq n$. The number $N$ of conjugates of $x$ in $S_n$ is given by the standard combinatorial formula 
$$N = \frac{n!}{(n-tq)! \, q^t \, t!} = \frac{1}{q^t} \binom{n}{t} \frac{(n-t)!}{(n-tq)!}.$$ The size of the conjugacy class $x^G$ in $A_n$ is either $N$ or $N/2$, depending on whether $C_{S_n}(x)$ contains odd permutations or not.

We proceed by analyzing first the case $q = 2$. Since $x$ is a product of $t$ disjoint transpositions with $t \geq 2$ even, there always exists an odd permutation centralizing $x$ (such as one of the transpositions forming $x$ itself). Thus, the class of $x$ in $S_n$ does not split in $A_n$, which implies 
$$|x^G| = N = \frac{n!}{(n-2t)! \, 2^t \, t!} = \frac{1}{2^t} \binom{n}{t} \frac{(n-t)!}{(n-2t)!}.$$
Direct verification for $n = 5, 6, 7$ (e.g. using GAP \cite{GAP4}) shows that no such conjugacy class size exists. For $n \geq 8$ and $t \geq 2$ even, the factor $n(n-1)\ldots (n-2t+1)$ of $N$ contributes at least $t(q-1)$ consecutive integers, so $N$ is always divisible by at least two distinct odd primes, which contradicts $|x^G| = 2p^a$.

Thus, for the remainder of the proof,  we assume that $q \geq 3$. As above, direct calculation confirms the theorem for $n=5,6,7$, so we assume $n\geq 8$.
 As $q$ is odd, any $q$-cycle is an even permutation. For $n \ge 8$, if $n = q$, the class splits and $|x^G| = N/2 = (n-1)!/2$, which is easily seen to be divisible by at least two distinct odd primes. The same happens when $n=q+1$ since $|x^G|=N/2=n(n-2)!/2$. Otherwise, if $n \geq q+2$, there always exists an odd permutation in $S_n$ centralizing $x$ (either by swapping two $q$-cycles if $t \ge 2$, or by transposing two fixed points if $t = 1$), so the class does not split and $|x^G| = N$.  We distinguish two possibilities according to the value of $t$.

When $t > 1$, the factor $(n-t)!/(n-tq)!$ contains the product of $t(q-1) \geq 4$ consecutive integers. Hence $N$ is divisible at least by two different odd prime numbers, so $|x^G| = N$ cannot equal $2p^a$. 
Finally, when $t = 1$, we have $$N = \frac{n(n-1)\cdots(n-q+1)}{q}.$$ If $q \geq 5$,  the numerator of $N$ is the product of $q$ consecutive integers. Since dividing by $q$ cancels at most the prime $q$ itself, this guarantees that $N$ is divisible by 8, a contradiction. If $q = 3$, then $$N = \frac{n(n-1)(n-2)}{3}.$$Now, if $n$ is even, then $n(n-2)$ is divisible by $8$, so $8$ divides $N$, which is impossible. Finally,  if $n \geq 9$ is odd, the product $n(n-1)(n-2)$ contains at least two distinct odd factors, $n$ and $n-2$, both greater than or equal to $7$; this guarantees that $N$ is divisible by at least two distinct odd primes. This completes the proof.
\end{proof}

\begin{theorem} \label{T2}
The sporadic simple groups do not possess any conjugacy class of size $2p^a$ for any prime $p$ and any integer $a\geq 0$.
\end{theorem}

\begin{proof} 
The theorem can be verified using the Atlas of Finite Groups  \cite{Con} or via the GAP Character Table Library \cite{GAP4} by computing the sizes of the conjugacy classes of the sporadic simple groups, which can be read off from their well-known character tables. Alternatively, the result can be deduced from  \cite[Theorem 1.1 and Table 3.2]{LL}. The only possible pairs $(G,H)$, where $G$ is a sporadic simple group and $H$ is a subgroup of $G$ of index $2p^n$, for some prime $p$ and some positive integer $n$, are $(M_{11}, A_6)$ and $(M_{22}, L_3(4))$. However, neither of these subgroups can be the equal to the centralizer of some element of $G$,  since both have trivial center.
\end{proof}

For the groups of Lie type, we base the proof of Theorem A on the classification of primitive groups of degree $p^a q^b$ obtained in~\cite{LL}. More precisely, we shall use Theorem~1.1 and Tables~4.1, 4.2, 4.3, 4.5, and~5.1, which provide the simple groups of Lie type having subgroups of index $p^a q^b$, where $p$ and $q$ are primes.  We shall use, following \cite{LL} the notation of the book \cite{KL}. Let $S$ be a Sylow $s$-subgroup of  a simple group of Lie type $G$ defined on a field of characteristic $s$, and for a positive integer $i$, let $L_i$ be the Levi subgroup.  Then $B=N_G(S)=SD$ is the Borel subgroup of $G$, and $P_i=BL_i=(O_s(P_i)L_i)D$ are parabolic subgroups with $O_s(P_1)L_i\lhd P_1$, where $D$ is the diagonal subgroup of $G$. The subgroup $L_1D$ is called the Levi complement (or Levi factor) of $P_i$.

\begin{theorem} \label{T3}
The simple groups of Lie type do not possess any conjugacy class of size $2p^n$ for any prime $p$ and any integer $n\geq 0$.
\end{theorem}

\begin{proof}
Let $G$ be a simple group of Lie type, and suppose there exists $x \in G$ such that the conjugacy class $x^G$ has size $|G :C_G(x)| = 2p^n$, where $p$ is a prime and $n \ge 1$ is an integer. Set $H = C_G(x)$. In view of the fact that we may assume $p\neq 2$,  the pair $(G, H)$ must appear in the classification given in ~\cite{LL}. More precisely, the candidates for $(G, H)$ are listed in Theorem~1.1 and Tables~4.1, 4.2, 4.3, 4.5, and~5.1 of~\cite{LL}. 

\medskip
We first observe that for all candidates in Table~4.2 through Table~5.1, whenever the subgroup $H$ has index $2p^n$, its center is trivial, contradicting the fact that $x\in Z(C_G(x))$ with $x\neq 1$. More specifically, the only candidate pairs in these tables whose indices  have the required form $2p^n$ are as follows: in   Table~4.2 the pairs $({\rm PSp}_4(3), 2^4: D_{10})$ and $({\rm Sp}_4(2), {\rm O}_4^-(2))$, with indices $2\cdot 3^4$ and $2\cdot 3$, respectively. However, ${\rm Sp}_4(2)\cong S_6$ is not simple, so we can exclude it.  In  Table~4.3, the unique pair is $({\rm U}_3(5), A_7)$ with index $2\cdot 5^2$. In Tables 4.4, 4.5 and 5.1 ,  a straightforward arithmetic check shows that no pairs yield an index of the form $2p^n$, because
 all of them are either odd numbers or divisible by $4$. Therefore, all cases appearing in Tables~4.2 through 5.1 are ruled out. 

\medskip
Consequently, $(G, C_G(x))$ must belong to the configurations listed in \cite[Table 4.1]{LL}, restricting our focus to $G \cong {\rm PSL}(m, q)$, with $q=s^f$ and $s$ prime. 
 A straightforward index verification rules out all entries in the upper section of that table, restricting our attention only to the families in the lower section, which we next discuss.

\medskip
(1) When $(G, C_G(x)) = ({\rm PSL}(2, 2^{2^i}), D_{2(2^{2^i} - 1)})$, an index computation gives $2^{2^i - 1} = 2$, which forces $i = 1$. Consequently, $G \cong {\rm PSL}(2, 4) \cong A_5$, a case previously eliminated in Theorem 2.1. 
Likewise, if $(G, C_G(x)) = ({\rm PSL}(2, 2^{2^i}), C_{2^{2^i} - 1})$, the index constraint forces $i = 0$, leading to $G \cong {\rm PSL}(2, 2) \cong S_3$, which contradicts the simplicity of $G$. By applying the same reasoning to the remaining families of parameters, we can readily eliminate $G \cong \mathrm{PSL}(2, 2^f)$, $\mathrm{PSL}(2, s^{2^i})$, $\mathrm{PSL}(2, s)$, and $\mathrm{PSL}(2, 2^{2^{i+1}})$, showing by elementary arithmetical arguments that the corresponding index of $H$ cannot equal $2p^n$, except in  the case $(\mathrm{PSL}(2, 5),  D_6)$, which has already been eliminated.

\medskip
We now examine the final four cases of Table 4.1 in  \cite{LL}, adopting its notation.

\medskip
(2)  Suppose that $(G, C_G(x)) = ({\rm PSL}(m, q), P_1)$, where $P_1$ is a maximal parabolic subgroup of index $\frac{q^m - 1}{q - 1} = 2p^n$ with $q = s^f$. According to \cite[Table 4.1]{LL}, either $m = r^2$ for a prime $r$ and $q=s$, or $m$ is prime and $f = m^i f_0$ (where $f_0 = 1$ or a prime). Notice that in both cases, the characteristic $s$ must be odd, as otherwise  $\frac{q^m - 1}{q - 1}$ would be odd, a contradiction.

Assume first that $m=2$.  In this case $P_1$ is a Borel subgroup of $G$ with structure $E_q\rtimes C_{(q-1)/2}$ (for $q$ being odd). Although this case is feasible from the  arithmetical point of view, when $q+1=2p^n$, the subgroup $P_1$  has trivial center, which immediately rules out this case since $x\in Z(P_1)$.

If $m = 4$. then the index becomes 
$$\frac{q^4 - 1}{q - 1} = (q + 1)(q^2 + 1) = (s^f + 1)(s^{2f} + 1) = 2p^n.$$  Since $s$ is an odd prime, both $s^f + 1$ and $s^{2f} + 1$ are even, implying $4 \mid (s^f + 1)(s^{2f} + 1) = 2p^n$,  a contradiction.

Therefore, we can assume $m > 2$ is an odd prime (or a square of an odd prime), and we know that $s$ is odd. Expanding the index yields
\[
\frac{q^m - 1}{q - 1} = \frac{s^{fm} - 1}{s^f - 1} = s^{f(m - 1)} + s^{f(m - 2)} + \dots + s^f + 1,
\] the right-hand side is a sum of $m$ odd integers, which is odd. This contradicts the fact that the index  is  an even integer. Therefore, this case is eliminated.

\medskip

(3) Suppose $(G, C_G(x)) = (\mathrm{PSL}(m, q), H)$, where $P_1$ is a maximal parabolic subgroup of $G$ with unipotent radical $O_s(P_1)$, and $H \le P_1$ satisfies $O_s(P_1)H = P_1$ and
\[
[G : H] = s^a \frac{q^m - 1}{q - 1} = 2p^n,
\]
with $q = s^f$ for a prime $s$. Note that 
\[
\frac{q^m - 1}{q - 1} = 1 + q + \dots + q^{m-1} \equiv 1 \pmod s,
\]
so $(q^m - 1)/(q - 1)$ is coprime to $s$. Since $p$ is an odd prime, the factor $2$ in $2p^n$ must be absorbed by $s^a$, forcing $s^a = 2$. Thus $s = 2$, $a = 1$, $q$ is even, $|G : P_1| = p^n$, and $|P_1 : H| = 2$.

We know that $P_1 = O_2(P_1) \rtimes L_1D$, where $L_1$ is the Levi subgroup. Moreover, if $m \ge 3$, then $P_1$ takes the form
\[
 P_1 = (q^{m-1} : \mathrm{GL}(m - 1, q)) / (Z(\mathrm{SL}(m, q)) \cap \mathrm{GL}(m - 1, q)),
\]
(see for instance \cite[Theorem 4.1.17]{KL}) and $L_1\cong {\rm PSL}(m-1,q)$ acts irreducibly on $O_2(P_1) \cong \mathbb{F}_q^{m-1}$. If $m = 2$, then $O_2(P_1) \cong \mathbb{F}_q^+$, and the abelian factor $L_1$ acts transitively on non-zero elements of $O_2(P_1)$.

Since $|P_1 : H| = 2$, we have $H \trianglelefteq P_1$, which implies that $H \cap O_2(P_1) \trianglelefteq P_1$. Then, $|O_2(P_1) : H \cap O_2(P_1)| = |O_2(P_1)H : H| = |P_1 : H| = 2$.
Thus, $H \cap O_2(P_1)$ is an $L_1$-invariant submodule (or subspace) of index 2 in $O_2(P_1)$. 

Now, if $m \ge 3$, then $|O_2(P_1)|=q^{m-1}\geq4$,
so $H\cap O_2(P_1)$ is non-trivial and proper,  contradicting the irreducibility of the action of $L_1$  on $O_2(P_1)$.
If \(m=2\), then $O_2(P_1)\cong \mathbb{F}_q^+$, and \(H\cap O_2(P_1)\) is an \(L_1\)-invariant subspace of index \(2\). Since \(L_1\) acts transitively on the non-zero elements of \(O_2(P_1)\), such an invariant subspace can exist only when \(q=2\). This yields $G\cong\mathrm{PSL}(2,2)$, which is not simple, also a contradiction.

\medskip

(4) Suppose $(G, C_G(x)) = (\mathrm{PSL}(m, q), O_s(P_1) L_1 R)$, where $q = s^f$ for a prime $s$, and $R<D$ with $|D:R|=r^b$ dividing $(q-1)$. The index formula yields 
\[
\frac{r^a (q^m - 1)}{q - 1} = 2p^n \quad \text{with } r^a \mid (q - 1),
\]
as indicated in  \cite[Table 4.1]{LL}. If $s = 2$, the condition $r^a \mid (q - 1)$ implies that $r$ is an odd prime. Moreover, given that $q$ is even, the quotient $(q^m - 1)/(q - 1)$ is odd. Thus, the left-hand side is a product of odd integers, contradicting the evenness of $2p^n$. Hence, $s$ is an odd prime.

Since $s \neq 2, p$ (as $s$ is odd and coprime to $\frac{q^m-1}{q-1}$, which is divisible by $p$)
we have that $C_G(x) = O_s(P_1) L_1 R$ contains a Sylow $s$-subgroup $S$ of $G$. In fact, we have $S\leq O_s(P_1) L_1$ since $O_s(P_1)L_1$ is a normal subgroup of $C_G(x)$ of index coprime to $s$. Hence, $x$ centralizes $S$. It is well known (see, for instance, \cite[Lemma 2.2]{BGL}) that $C_G(S) = Z(S)$, which implies $x \in Z(S)$. Furthermore, we have  $Z(S) \subseteq O_s(P_1)$ \cite[Lemma 3.3]{BGL}. For the family  ${\rm PSL}(m,q)$, this inclusion can be verified easily. Since $S \le P_1$,  without loss of generality,  $S$ can be assumed to be a subgroup of upper unitriangular matrices modulo $Z({\rm SL}((m,q)$. Then $Z(S)$ consists of matrices with a single non-zero off-diagonal entry at position $(1, m)$, which is clearly contained in the unipotent radical $O_s(P_1)$ (comprising matrices with arbitrary entries along the first row).

Now, as mentioned in case (3), if $m \ge 3$ then $L_1$ acts irreducibly on $O_s(P_1) \cong \mathbb{F}_q^{m-1}$, while if $m = 2$ then $L_1$ acts transitively on the non-zero elements of $O_s(P_1) \cong \mathbb{F}_q^+$. In both cases, $L_1$ acts on $O_s(P_1)$ without non-trivial fixed points, meaning that $C_{O_s(P_1)}(L_1)=1$. Since $x \in O_s(P_1)$ is non-trivial and $L_1 \le C_G(x)$, we conclude that this case cannot occur.

\medskip
(5) Finally, if $(G, \mathbf{C}_G(x)) = ({\rm PSL}(4, 2^{2^i}), P_1)$, the index of $P_1$ is given by
\[
\frac{(2^{2^i})^4 - 1}{2^{2^i} - 1} = (2^{2^i} + 1)(2^{2^{i+1}} + 1).
\]
Since $2^{2^i} + 1$ and $2^{2^{i+1}} + 1$ are both odd integers for any $i\geq 1$, their product is strictly odd, contradicting the requirement that the index equals $2p^n$. This case completes the proof.
\end{proof}

\begin{remark}\label{arithmetic_viability}
The proof of Theorem 2.3 shows that, at the arithmetical level,  the condition
$|G:H|=2p^n$
 eliminates most of the configurations in Table 4.1 through Table 5.1 in \cite{LL} without imposing the additional requirement that $H=C_G(x)$. After these eliminations, the only cases that remain arithmetically viable are the following:
\begin{enumerate}
\item[(a)] \(({\rm PSp}_4(3),\,2^4:D_{10})\) and \(({\rm U}_3(5),\,A_7)\). Observe that the case \(({\rm PSL}(2,5),S_3)\), which is isomorphic to \((A_5,S_3)\), has already been considered in Theorem 2.1).

\item[(b)] \(({\rm PSL}(2,q),P_1)\), where \(q=s^f\) is odd and
$
\frac{q^2-1}{q-1}=q+1=2p^n,
$
with \(P_1\) a parabolic subgroup. Note that $s\neq p$.

\item[(c)] \(({\rm PSL}(m,q),O_s(P_1)L_1R)\), where \(q=s^f\) is odd and
$r^a\frac{q^m-1}{q-1}=2p^n,$
with $r$ prime such that \(r^a\mid(q-1)\), and \(P_1\) a parabolic subgroup with unipotent radical \(O_s(P_1)\), Levi factor \(L_1\), and $R<D$ where $|D:R|=r^a$. Note that $s\neq p$.
\end{enumerate}
\end{remark}

\section{Proofs of Theorems B and C}

We start by proving Theorem C.  Our proof is also based on the classification of primitive groups of degree $p^a q^b$ obtained in \cite{LL}. Specifically, we utilize Theorem 1.1 and Tables 3.1, 3.2 , 4.1, 4.2, 4.3, 4.5, and 5.1 of \cite{LL}, which list the simple groups, alternating, sporadic, and of Lie type, containing subgroups $H$ of index $p^a q^b$. We shall also appeal to the structure of the automorphism groups of finite simple groups, particularly that of ${\rm PSL}(m, q)$, for which we  shall use \cite{BG}.

We also base our proof on Kazarin's $p^a$-Theorem and the following consequence of it. 

\begin{lemma} \label{4.1} Let $G$ be a finite group and $\alpha$ one of its automorphisms. If $C_G(\alpha)$ contains a Sylow $r$-subgroup of $G$ for all primes $r \in \pi(G) \setminus \{p\}$, then $\alpha$ induces the identity automorphism on $G / S(G)$ (where $S(G)$ denotes the solvable radical of $G$). In particular, a non-abelian simple group $G$ has no non-trivial automorphism $\alpha$ such that $|G:C_G(\alpha)|$ is a prime power.
\end{lemma}
\begin{proof}
    The first statement of the Lemma corresponds to  \cite[Corollary 1]{Kaz}, from which the second assertion follows directly.
\end{proof}

 \begin{theorem}[Theorem C]
     Let $G$ be a non-abelian finite simple group and suppose there exists $a\in {\rm Aut}(G)$ such that $|G:C_G(a)|=2p^n$ with $p$ prime and $n\geq 1$. Then $G\cong A_5$, $p=5 $, $n=1$,  and $a$ induces the unique outer automorphism of $A_5$ (corresponding to a transposition of $S_5$).   
    \end{theorem}
\begin{proof}

We first show that, without loss of generality, we may assume that $a$ has prime order. Suppose that $o(a) = k > 1$ is not a prime number, and let $r$ be a prime factor of $k$. Consider $b = a^{k/r} \in {\rm Aut}(G)$, which has prime order $r$. We have $C_G(a) \le  C_G(b) \le G$. Consequently, $|G :C_G(b)|$ divides $|G : C_G(a)| = 2p^n$, so $|G :C_G(b)| = p^m$  or$|G :C_G(b)| = 2p^m$ for some $0 \le m \le n$. 
Now, if $|G : C_G(b)| = p^m$ or $2$, Lemma \ref{4.1} implies that a non-trivial automorphism  cannot fix a subgroup of prime-power index in a non-abelian simple group, forcing $b = 1$, a contradiction. Thus, $|G :C_G(b)| = 2p^m$ with $m\geq 1$. Henceforth, we may assume that $a$  has prime order. 

\medskip
 We can assume further that $a\not \in$ Inn $(G)$. To the contrary, $a$ corresponds to conjugation by some element $g \in G$, which implies that $C_G(a) = C_G(g)$. Consequently,  $g^G$ has size $2p^n$, directly contradicting Theorem~A.
We divide the proof into three cases attending to the CFSG. Our strategy consists in checking each candidate pair $(G, H)$ in \cite[Theorem 1.1]{LL}  and verifying that, apart from the  exceptional case in the statement, $H$ cannot occur as the fixed-point subgroup of any automorphism of prime order of $G$.
\medskip

(1) Suppose that $G \cong A_m$ with $m \ge 5$. Suppose first that $m\leq 8$. By inspecting \cite[Table 3.1]{LL}, we obtain exactly four candidate pairs for $(G, C_G(a))$: the pairs $(A_5, D_{10})$, $(A_6, 3^2:4)$ and $(A_6, A_5)$ from the upper part of the table; and the pair $(A_5, S_3)$ from the lower part. We analyze these possibilities case by case.   
If $G \cong A_5$, then   $a$ induces the sole outer automorphism of $A_5$ of order~$2$, corresponding to a transposition in $S_5 \cong G \rtimes \langle a \rangle$. If $C_G(a) \cong D_{10}$,  then the centralizer of $a$ in the semidirect product $G \rtimes \langle a \rangle\cong S_5$ would contain a subgroup isomorphic to $C_G(a) \times \langle a \rangle \cong D_{10} \times C_2$, which is impossible.  If $C_G(a) \cong S_{3}$,  then $|S_5: C_{S_5}(a)|=|G : C_G(a)| =2 \cdot 5$, which is of the form $2p^n$ with $p = 5$ and $n = 1$. This case actually occurs, since  $C_{S_5}(a)\cong S_3\times C_2$, yielding the configuration stated in the theorem. Consider now the case $G \cong A_6$, then $\mathrm{Out}(A_6) \cong C_2 \times C_2$. The group $A_6$ possesses exactly three outer automorphisms of order $2$, corresponding to the three extensions of index $2$ in $\mathrm{Aut}(A_6)$; namely, $S_6$, $M_{10}$, and $\mathrm{PGL}(2,9)$. However, none of these groups possesses a subgroup isomorphic to $C_G(a) \times \langle a \rangle \cong (3^2:4) \times C_2$ or $A_5\times C_2$ (readily checked using \cite{GAP4}), which rules out the  pairs $(A_6, 3^2:4)$ and $(A_6, A_5)$.

Now suppose that $m > 8$, so we are in the lower part of \cite[Table 3.1]{LL}. We observe that the only case in which $A_m$ possesses a subgroup $H$ of index $2p^n$ is when $m = 2p^n$ and $H = C_G(a) \cong A_{m-1}$. In this setting, since $\mathrm{Out}(A_m) \cong C_2$ for $m > 8$, and $a$  has prime order, $a$ is induced by  a transposition in $S_m$. However, the fixed-point subgroup of a transposition of $S_m$  in $A_m$ is isomorphic to $A_{m-2} \times C_2\cong S_{m-2}$, which contradicts $C_G(a) \cong A_{m-1}$. Thus, no further solutions arise.

 \medskip
(2) Consider now the case where $G$ is a sporadic simple group. According to \cite[Table 3.2]{LL}, the only candidate pairs for $(G, C_G(a))$ are $(M_{11}, A_6)$ and $(M_{22}, \mathrm{PSL}(3,4))$.  
In view of $\mathrm{Out}(M_{11}) = 1$, the group $M_{11}$ only has inner automorphisms, so no solution arises from this case.  Thus, we are left with $G \cong M_{22}$. Since $\mathrm{Out}(M_{22}) \cong C_2$, then $a$ induces the sole outer automorphism of $M_{22}$. If $C_G(a) \cong \mathrm{PSL}(3,4)$, then the centralizer of $a$ in the semidirect product $G \rtimes \langle a \rangle$ would contain a subgroup isomorphic to $C_G(a) \times \langle a \rangle \cong \mathrm{PSL}(3,4) \times C_2$. However, consulting the maximal subgroups of $\mathrm{Aut}(M_{22}) \cong M_{22}:2$ in \cite{Con}, no such subgroup exists, allowing us to conclude that $G\not \cong M_{22}$.

\medskip

(3) Suppose that $G$ is a simple group of Lie type.   We follow the same strategy as in the proof of Theorem 2.3, relying on the classification of the simple groups having  subgroups of index $p^a q^b$ given in \cite[Tables 4.1--4.5 and 5.1]{LL}. In view of Remark \ref{arithmetic_viability}, and following the same notation, the only cases that remain arithmetically viable are:
\begin{itemize}
\item[\rm (3.1)]   $({\rm PSp}_4(3), 2^4: D_{10})$  and $({\rm U}_3(5), A_7)$. 
    \item[\rm (3.2)] $(\mathrm{PSL}(2, q), P_1)$, with $q=s^f$ odd, and $(q^2-1)/(q-1)=q+1 =2p^n$; and
    \item[\rm (3.3)] $({\rm PSL}(m,q), {\bf O}_s(P_1)L_1R)$, where $r^b\frac{q^m-1}{q-1}=2p^n$ and $q=s^f$ odd. 
\end{itemize}

(3.1) Suppose first $(G, C_G(a)) = (\mathrm{PSp}_4(3), 2^4 : D_{10})$. 
    Here, $\mathrm{Out}(G) \cong C_2$, so $G$ possesses a sole extension of degree $2$, namely $\mathrm{PGSp}_4(3) \cong G \rtimes \langle a \rangle$. If there were  such an element $a \in \mathrm{Aut}(G) \setminus \mathrm{Inn}(G)$ of order $2$, then its centralizer in $\mathrm{PGSp}_4(3)$ would contain a subgroup isomorphic to $C_G(a) \times \langle a \rangle\cong 2^4 : D_{10}\times 2$.  However, inspecting  in \cite{Con} (or using GAP \cite{GAP4}), we find that no such subgroup exists. Hence, this pair is ruled out.

    Consider now the case $(G, C_G(a)) = (\mathrm{U}_3(5), A_7)$. 
    Then $\mathrm{Out}(\mathrm{U}_3(5)) \cong S_3$ and contains a unique conjugacy class of outer involutions corresponding to graph automorphisms, and a class of outer diagonal automorphisms of order $3$. Assume first $o(a)=2$ and consider  $G \rtimes \langle a \rangle \cong \mathrm{U}_3(5) \rtimes C_2$.  Since $a$ centralizes $C_G(a)\cong A_7$, then the extension $G \rtimes \langle a \rangle$ must contain $A_7 \times C_2$. However, the unique degree-$2$ extension of $A_7$ embedded in $\mathrm{Aut}(\mathrm{U}_3(5))$ extends $A_7$ to $S_7$, contradicting $a$ acting trivially on $A_7$.   Similarly, if $a$ induces an outer automorphism of order $3$, then $G \rtimes \langle a \rangle \cong \mathrm{PGU}_3(5)$. If $a$ centralizes $A_7$, then $\mathrm{PGU}_3(5)$ would contain a subgroup isomorphic to $A_7 \times C_3$. However, the centralizer of $A_7$ in $\mathrm{PGU}_3(5)$ is trivial. Thus,  this case cannot occur either.

\medskip
(3.2) and (3.3) In both cases, $s$ is odd. Observe that $C_G(a)$ contains a Sylow $s$-subgroup of $G$, since $s\neq 2$ and $s\neq p$. Let $T = \mathrm{PGL}(m, q)$, so we have $G \le T$. Since $a$ is an automorphism of prime order of $G$, we can determine the order of $C_T(a)$ using \cite[Proposition 3.2.16 and Table B.3]{BG}. Specifically, it suffices to compare the $s$-part $|C_T(a)|_s$ with the $s$-part $|G|_s=q^{m(m-1)/2}$, demonstrating that $|C_T(a)|_s < |G|_s$, which contradicts that $C_G(a)$ contains a Sylow $s$-subgroup of $G$. 

We distinguish the various types of automorphisms  of prime order following the classification given in [2, Table B.3];  namely, $\texttt{u}$ (inner), $\texttt{s}$ (diagonal),  $\texttt{f}$  (field), $\texttt{gf}$ (graph-field) and $\texttt{g}$ (graph).  Type ${\texttt u}$ is excluded because $a$ is non-inner. If $a$ is of type $\texttt{s}$ with the condition $i > 1$, then $a$ is a diagonal automorphism, and the centralizer order is given by
$$|C_T(a)| = \frac{1}{q-1} |\mathrm{GL}(e, q)| \prod_j |\mathrm{GL}(a_j, q^i)|$$
(where $e$ and $a_j$ are defined as in \cite[Table B.3]{BG}). Using standard formulas for the orders of general linear groups, we have
$$|\mathrm{GL}(e, q)| = q^{\frac{e(e-1)}{2}} \prod_{t=1}^e (q^t - 1)\,\, {\rm and} \, \,|\mathrm{GL}(a_j, q^i)| = q^{\frac{a_j(a_j-1)i}{2}} \prod_{t=1}^{a_j} (q^{it} - 1).$$
Now, since $m = e + i \sum_j a_j$, the $s$-part of $|C_T(a)|$ evaluates to
$$|C_T(a)|_s = q^{\frac{e(e-1) + i \sum_j a_j(a_j-1)}{2}}.$$
Given that $i > 1$, we obtain
$$e(e - 1) + i \sum_j a_j(a_j - 1) < \left(e + i \sum_j a_j\right)^2 - \left(e + i \sum_j a_j\right) = m^2 - m,$$
which guarantees that $|C_T(a)|_s < q^{\frac{m(m-1)}{2}} = |G|_s$. Analogous order comparisons apply to rule out the remaining configurations of automorphisms of type $\texttt{s}$ and are omitted.

\medskip
The other types of automorphisms: $\texttt{f}$, $\texttt{gf}$ and $\texttt{g}$, can be handled similarly. For instance, consider the case where $a$ is of type $\texttt{g}$ with $m\geq 4$ even, $o(a)=2$ and $C_T(a) = \mathrm{PGSp}_m(q)$. Then we have
$$|\mathrm{PGSp}_m(q)| = q^{(m/2)^2} \prod_{i=1}^{m/2} (q^{2i} - 1),$$
so its $s$-part is given by $|C_T(a)|_s = q^{m^2/4}$. Since $m \ge 4$, it holds that $m^2/4 < m(m-1)/2$, which again implies $|C_T(a)|_s < |G|_s$. The field and graph-field automorphisms of prime order $r$ satisfy $|C_T(a)|_s \le q^{m(m-1)/(2r)}$, which is likewise strictly smaller than $\vert{}G\vert{}_s$.

Thus, in all cases for $a$, we conclude that $|C_T(a)|_s < |G|_s$. Since $|C_G(a)|_s \le |C_T(a)|_s$, it follows that $|C_G(a)|_s < |G|_s$, completing the contradiction.
\end{proof}

For proving Theorem B, we first present the following preliminary result of elementary character.

\begin{lemma} \label{cent}
Suppose that a finite group $A$ acts on a finite group $G$ which allows
a direct decomposition $G = H_1\times \cdots \times H_n$, that is invariant under $A$, i.e.,
$H^a_i\in \{H_1,...,H_n\}$ for all $a \in A$ and all $i \in \{1,...,n\}$. Assume further that
A acts transitively on $\{H_1,...,H_n\}$. Let $H \in  \{H_1,...,H_n\}$ and $B = {\bf N}_A(H)$. Then ${\bf C}_G(A)\cong {\bf C}_H(B)$.
\end{lemma}

\begin{proof}
This is \cite[8.1.6]{KS}.
\end{proof}

\begin{theorem}[Theorem B] Let $G$ be a finite group and suppose that $1\neq g\in G$ is an element of odd order with $|g^G|=2p^n$, where $p$ is a prime and $n\geq 0$. Then $g^G$ generates a solvable normal subgroup of $G$.
\end{theorem}

\begin{proof}

 We first note that if $p = 2$ or $n = 0$, then $|g^G|$ is a power of $2$, and Kazarin's Theorem applies directly. Hence, assume $p$ is odd and $n \ge 1$. Let $G$ be a counterexample of minimal order, and choose $g\in G$ of minimal order among all counterexamples in $G$.

\medskip
We make the following observation: if $1\neq N$ is a normal subgroup of $G$ such that $g \notin N$, the subgroup $\langle (gN)^{G/N} \rangle$ is solvable in $G/N$. Indeed, the image $gN$ has (odd) order dividing $o(g)$, and its conjugacy class size $|(gN)^{G/N}|$ divides $|g^G| = 2p^n$. If $N \neq 1$, then $|G/N| < |G|$, so the minimality of $|G|$ (together with Kazarin's Theorem if the class size is a prime power) implies that $\langle (gN)^{G/N} \rangle = \langle g^G \rangle N / N$ is solvable.

\medskip
We first claim that $o(g)$ is prime. Suppose, on the contrary, that $o(g)$ is composite. Choose $q \mid o(g)$ and set $1 \neq y = g^{o(g)/q}$. Then $o(y) = q < o(g)$ and $C_G(g) \le C_G(y)$, so $|y^G|$ divides $|g^G| = 2p^n$. By Kazarin's Theorem and the minimality of $o(g)$, we have that $1\neq N = \langle y^G \rangle$ is solvable. Since $g^{o(g)/q} = y \in N$, we must have $g \notin N$, otherwise $\langle g^G\rangle$ would be solvable. By our initial observation applied to $N$, the factor group $\langle g^G \rangle N / N$ is solvable. Since $N$ is solvable, it follows that $\langle g^G \rangle$ is solvable, a contradiction. Thus $o(g)$ is prime.

\medskip
Put $M = \langle g^G \rangle$. We next prove that $M = G$. Suppose, on the contrary, that $M < G$. Since $M\lhd G$, then $|g^M|$ divides $|g^G| = 2p^n$, and Kazarin's Theorem together with the minimality of $|G|$ imply that $\langle g^M \rangle$ is solvable. This means that $g$ lies in the solvable radical of $M$, and consequently, in the solvable radical of $G$. This forces $\langle g^G \rangle$ to be solvable, a contradiction. This proves $M = G$.

\medskip
Now, by Theorem A, we know that $G$ is not simple. We claim that every minimal normal subgroup of $G$ is non-solvable. Suppose $N$ is a solvable minimal normal subgroup of $G$. Since $M = G$, we have $g \notin N$. Applying our initial observation to $N$, we obtain that $G/N = \langle g^G \rangle N / N$ is solvable. The solvability of $N$ then implies that $G$ is solvable, a contradiction. This proves that every minimal normal subgroup of $G$ is non-solvable.

\medskip
Thus, we can let $N<G$ be a minimal normal subgroup of $G$, which can be written as
$N=S_1\times\cdots\times S_t$,
where the $S_i$ are isomorphic non-abelian simple groups. Since
$M=G$, we have $g\notin N$. Set
$K=N\langle g\rangle$. It is clear that
$K=NC_K(g)$. 
Consequently,
\[
|g^K|
=
|K:C_K(g)|
=
|N:C_N(g)|
=
|NC_G(g):C_G(g)|,
\]
and thus $|g^K|$ divides $|g^G|=2p^n$. Next, we  prove that $K=G$.
Suppose to the contrary that $K<G$. By the minimality of $|G|$, together with
Kazarin's Theorem when $|g^K|$ is a prime power, the subgroup
$L=\langle g^K\rangle$
is solvable. Moreover, $L\lhd K$. Since $N\lhd K$, the subgroup
$L\cap N$ is solvable and normal in $N$, so it follows that $L\cap N=1$. This means that $g\in C_G(N)$. Since $C_G(N)\lhd G$, every conjugate of $g$ also
centralizes $N$, so the fact that $G=\langle g^G\rangle$ yields
$G=C_G(N)$, that is, $N\leq Z(G)$, a contradiction. Hence,
$K=G$, or equivalently $G=N\langle g\rangle$.

\medskip
We now prove that $N$ is simple. Since $N$ is a minimal normal subgroup
of $G=N\langle g\rangle$, the cyclic group  $\langle g\rangle$ acts
transitively by conjugation on the set
$\{S_1,\ldots,S_t\}.
$
Set $S=S_1$, and note that if $t>1$, then
the stabilizer $N_{\langle g\rangle}(S)=1$ since $g$ has prime order. 
By Lemma \ref{cent},
$C_N(g)=C_N(\langle g\rangle)\cong C_S(1)=S$.
Therefore,
\[
2p^n
=|G:C_G(g)|
=|N:C_N(g)|
=|S|^{t-1}.
\]
 This is impossible,
because $S$ is a non-abelian simple group, so $|S|$ is divisible by at least three distinct primes. Thus $t=1$, meaning that $N$ is simple.

\medskip
Finally, conjugation by $g$ induces an automorphism
$a\in\operatorname{Aut}(N)$. Since
$C_N(a)=C_N(g),
$
we have
\[
|N:C_N(a)|
=
|G:C_G(g)|
=
2p^n.
\]
Hence Theorem~C applies. In particular, it follows that
 $a$ has even order. However, since $a$ is induced by $g$, $o(a)$ must divide $o(g)$, which is impossible since $o(g)$ is odd. This final contradiction proves the theorem.
\end{proof}

\section{The non-solvable case}
In the preceding section, we established that elements of odd order in a group $G$ with conjugacy class size $2p^n$ necessarily lie in the solvable radical of $G$. When an element $x\in G$ has even order, however, $\langle x^G \rangle$ need no longer be solvable, as illustrated by $S_5$ in Theorem C. The objective of this section is to determine the structure of the subgroup generated by $x^G$, when it is non-solvable.

\begin{theorem}[Theorem D]
Let $G$ be a finite group containing a conjugacy class $g^G$ of size $2p^n$, where $p$ is prime and $n\geq 1$,  such that $\langle g^G\rangle$ is non-solvable. Then $p=5$, and 
$\langle g^G \rangle/S(\langle g^G \rangle)$ is isomorphic to $S_5$.
\end{theorem}

\begin{proof}
 Let $S(G)$ denote the solvable radical of $G$, and set $\overline{G} = G/S(G)$ and $\overline{N} = \langle g^G \rangle S(G)/S(G)$. Since $\langle g^G\rangle $ is non-solvable then $\overline{N}\neq 1$. Let $\overline{L} = \text{Soc}(\overline{G}) = L_1 \times \cdots \times L_k$ be the socle of $\overline{G}$, where each $L_i$ for $i=1,\ldots ,k$ is a minimal normal subgroup of $\overline{G}$. Since $S(\overline{G}) = 1$, each $L_i$ is a direct product of non-abelian simple groups.

\medskip
We first prove that $\overline{N}$ centralizes all factors of $\overline{L}$ except one. It is clear that $C_{\overline{G}}(\overline{L}) = 1$, so $\overline{N}$ cannot centralize all $L_i$. Assume, for instance, that $L_1$ is not centralized by $\overline{N}$. This implies that $\overline{g}$ does not centralize $L_1$ either; otherwise, all conjugates of $\overline{g}$ in $\overline{G}$ would centralize $L_1$, and thus $\overline{N}$ as well. On the other  hand, as $\overline{g}$ normalizes every $L_i$, we have
$$1 \neq |\overline{L} : C_{\overline{L}}(\overline{g})| = \prod_{i=1}^k |L_i : C_{L_i}(\overline{g})|,$$
and note that this index divides $|\overline{G}:C_{\overline{G}}(\overline{g})|$, which divides $|g^G| = 2p^n$. Therefore, for every $i$, $|L_i : C_{L_i}(\overline{g})|$ must divide $2p^n$. But these indices cannot be prime powers by Lemma \ref{4.1}, so all of them are equal to $1$ except one, namely $|L_1 : C_{L_1}(\overline{g})|$, which must be equal to $2p^a$ for some $1\leq a \leq n$. As a result, $\overline{g}$ centralizes $L_i$ for all $i = 2, \dots, k$, as required.

\medskip
Write $L_1=S_1\times \ldots \times S_r$, where the $S_i$ are isomorphic non-abelian simple groups. Next we examine the action of $\overline{g}$ on $\{S_1, \dots, S_r\}$. We prove that this action, which is not necessarily transitive, cannot permute them with any orbit of length $t \geq 2$. Suppose, for instance, that $\{S_1, \dots, S_t\}$ is the orbit of $S_1$ under the action of $\overline{g}$, and set $T=S_1\times\ldots\times S_t$. By  applying Lemma 3.2 to the action of $\overline{g}$ on $T$, we obtain $C_{T}(\overline{g}) \cong C_{S_1}(N_{\langle \overline{g} \rangle}(S_1)) = C_{S_1}( \overline{g} ^s )$ for some $s\geq 1$. Then
$$1\neq |T : C_T(\overline{g})| =\frac{|S_1\times \ldots \times S_t|}{|C_{S_1}(\overline{g}^s)|}= |S_1 : C_{S_1}(\overline{g}^s)| \prod_{j=2}^t |S_j|.$$
Moreover, because $T$ is subnormal in $\overline{G}$, this index divides $|\overline{G}:C_{\overline{G}}(\overline{g})|$, and in turn $ |G:C_G(g)|$. But it is divisible by $|S_1|^{t-1}$. Thus, if $t>1$, this factor introduces at least three distinct prime divisors in $|G:C_G(g)|=2p^n$, which is impossible. Therefore, $t = 1$, which means that  $\overline{g}$ normalizes $S_i$ for all $i=1,\ldots,r$. We then have 
$$1\neq |{L_1} : C_{L_1}(\overline{g})| = \prod_j |S_j : C_{S_j}(\overline{g})|.$$ If $\overline{g}$ acted non-trivially on two or more simple factors, say $S_1$ and $S_2$, the product of the corresponding centralizer indices, which cannot be prime powers by Lemma \ref{4.1},  would have more than two even factors, leading to $4\mid |g^G|$, a contradiction. Hence, $\overline{g}$ acts non-trivially on exactly one simple factor, say $S_1$,  and trivially on $S_j$ for $2\leq j \leq r$.
\medskip

We have already shown that $\overline{N}$ centralizes all $L_i$ for $2\leq i \leq k $ and all  $S_j$ for $2\leq j \leq r$. Consequently, $\overline{N} \cap L_j \le Z(L_j) = 1$ for $2\leq j \leq k$. It follows that $L_1$ is the only factor of $\overline{L}$ that can be contained in $\overline{N}$; in fact, as $L_1$ is minimal normal in $\overline{G}$, we obtain $\overline{N}\cap \overline{L}=L_1$. Now, if $r>1$, it would follow that $L_1$ centralizes $S_j$ for all $2\leq j\leq r$, which is impossible. Therefore, $r=1$, so $L_1=S_1$. In particular, $S_1$ is normal in $ \overline{N}$. Moreover, we have 
 $C_{\overline{N}}(S_1) \leq C_{\overline{N}}(\overline{L}) \leq C_{\overline{G}}(\overline{L}) = 1$, and hence $\overline{N}$ embeds into
 $\text{Aut}(S_1)$. Thus,
$$S_1 \leq \overline{N} \leq \text{Aut}(S_1).$$
Furthermore, since $S_1$ is subnormal in $\overline{G}$, then $|S_1 : C_{S_1}(\overline{g})|$ divides $|\overline{G} : C_{\overline{G}}(\overline{g})|$, which in turn divides $|g^G| = 2p^n$.  By Lemma \ref{4.1}, $\vert{}S_1 : C_{S_1}(\overline{g})\vert{}$ cannot be a prime power, and as a consequence,  $|S_1 : C_{S_1}(\overline{g})|=2p^b$ with $1\leq b\leq n$. Applying Theorem C to the non-trivial action of $\overline{g}$ on $S_1$,  we conclude that $S_1 \cong A_5$ and $p = 5$. Therefore, $A_5 \leq \overline{N} \leq \text{Aut}(A_5) \cong S_5$, so  $\overline{N}$ must be isomorphic to either $A_5$ or $S_5$. However, since the action of $\overline{g}$ on $A_5$ induces an outer automorphism, $\overline{N}$ cannot be contained in $\text{Inn}(A_5) \cong A_5$. We conclude that 
$$\frac{\langle g^G \rangle}{S(\langle g^G \rangle)} \cong \frac{\langle g^G \rangle S(G)}{S(G)} \cong S_5.$$
This completes the proof. 
\end{proof}

\begin{example} It is worth pointing out that the exponent $n$ in the conjugacy class size $|g^G| = 2 \cdot 5^n$ in Theorem D need not be restricted to $n = 1$. Indeed, for any integer $n \ge 1$, one can construct a group with a conjugacy class of size $2 \cdot 5^n$ by taking a direct product with a suitable solvable group $H$. Specifically, let $H$ be a solvable group containing an element $z \in H$ (of any order), such that its conjugacy class size in $H$ is a power of 5, say $|z^H| = |H : C_H(z)| = 5^{n-1}$ (an explicit construction is omitted, as the existence of such groups is guaranteed). Taking $G = S_5 \times H$ and defining $x = (\tau, z) \in G$, where $\tau \in S_5$ is a transposition, we have $C_G(x) = C_{S_5}(\tau) \times C_H(z)$. Since $|\tau^{S_5}| = 10 = 2 \cdot 5$, the size of the conjugacy class of $x$ in $G$ is given by
\[
|x^G| = |G : C_G(x)| = |S_5 : C_{S_5}(\tau)| \cdot |H : C_H(z)| = (2 \cdot 5) \cdot 5^{n-1} = 2 \cdot 5^n.
\]
Notice that $o(x)=\operatorname{lcm}(2, o(z))$, so it is even.  In addition, since $H$ is solvable, the non-solvable normal subgroup generated by $x^G$, namely $\langle x^G \rangle = \langle \tau^{S_5} \rangle \times \langle z^H \rangle = (A_5\rtimes \langle \tau\rangle) \times \langle z^H \rangle$, satisfies that $\langle x^G\rangle/S(\langle x^G\rangle)$ is isomorphic to $S_5$, in full agreement with Theorem D.
\end{example}

\begin{example}
The following example demonstrates that in Theorem~D, the subgroup 
$S(\langle g^G \rangle)$ need not be a direct factor of $\langle g^G \rangle$.

Consider the standard wreath product  $G=C_5 \wr S_5$, and 
write $G = V \rtimes S_5$, with $V=C_5^5$ being the corresponding base group. 
Let $\tau = (1\,2) \in S_5$. 
The centralizer of $\tau$ in $G$ is given by $C_G(\tau) = C_V(\tau) \rtimes C_{S_5}(\tau)$. 
We have $C_{S_5}(\tau) = \langle (1\,2) \rangle \times S_3$, and  $|C_V(\tau)| = 5^4$. 
Therefore, $|C_G(\tau)| = 5^4 \cdot 12$, which yields $|\tau^G| =|G:C_G(\tau)|= 50=2\cdot 5^2$. Furthermore, it is straightforward to see that $\langle \tau^G \rangle = G$ 
and $S(G)=V$, which is not a 
direct factor of $G$.
\end{example}

\bigskip
\noindent
{\bf Acknowledgements}
This work was supported by the National Nature Science Fund of China (No.  12071181). The second author would like to express his sincere gratitude for the warm hospitality received during his research stay at the Universitat Jaume I, Spain, in autumn 2026.

\bibliographystyle{plain}

\begin{thebibliography}{1}


\bibitem{BGL}
N. Burgoyne, R. Griess, and R. Lyons,
\emph{Maximal subgroups and automorphisms of Chevalley groups},
Pacific J. Math. \textbf{71} (1977), no. 2, 365--403.

\bibitem{BG}
T. C. Burness and M. Giudici,
\emph{Classical Groups, Derangements and Primes},
Australian Mathematical Society Lecture Series, vol. 25,
Cambridge University Press, Cambridge, 2016.



\bibitem{Con}
J. H. Conway, R. T. Curtis, S. P. Norton, R. A. Parker, and R. A. Wilson,
\emph{Atlas of Finite Groups},
Oxford University Press, London, 1985.

\bibitem{Huppert}
B. Huppert, \textit{Character Theory of Finite Groups}, De Gruyter Expositions in Mathematics, vol. 25, Walter de Gruyter \& Co., Berlin, 1998.

\bibitem{Kaz}
L. S. Kazarin,
\emph{On groups with a conjugacy class of prime power order},
Mat. Sb. (N.S.) \textbf{115(157)} (1981), no. 1(5), 111--125;
English translation: Math. USSR-Sb. \textbf{43} (1982), no. 1, 99--111.

\bibitem{KL}
P. B. Kleidman and M. W. Liebeck,
\emph{The Subgroup Structure of the Finite Classical Groups},
LMS Lecture Note Series, vol. 129,
Cambridge University Press, Cambridge, 1990.


\bibitem{KS}
H. Kurzweil and B. Stellmacher,
\emph{The Theory of Finite Groups. An Introduction},
Springer-Verlag, Berlin-Heidelberg-New York, 2004.


\bibitem{LL}
C. H. Li and X. H. Li,
\emph{On permutation groups of degree a product of two prime powers},
Comm. Algebra \textbf{42} (2014), no. 11, 4722--4743.

\bibitem{GAP4}
The GAP Group, 
\emph{GAP -- Groups, Algorithms, and Programming, Version 4.13.1}; 
2024, https://www.gap-system.org.
\end{thebibliography}

\end{document}